\documentclass[11pt]{amsart}
\usepackage[T1]{fontenc}
\usepackage{amsmath,amssymb,amsthm}
\usepackage{etoolbox}
\makeatletter
\patchcmd{\@setauthors}{\MakeUppercase{\authors}}{\authors}{}{%
  \PackageError{authorcase}{Could not patch amsart author capitalization}{}%
}
\makeatother
\usepackage[margin=1.1in]{geometry}
\usepackage[colorlinks=true,linkcolor=blue,citecolor=blue,urlcolor=blue]{hyperref}
\hypersetup{
  pdftitle={On the Laugesen--Siudeja Conjecture: A Sharp Lower Bound and Failure of the Upper Bound},
  pdfauthor={Zhengchao Ji}
}
\numberwithin{equation}{section}
\newtheorem{theorem}{Theorem}[section]
\newtheorem{lemma}[theorem]{Lemma}
\newtheorem{proposition}[theorem]{Proposition}
\newtheorem*{conjecture}{Conjecture}
\newtheorem*{openproblem}{Open problem}
\theoremstyle{remark}
\newtheorem{remark}[theorem]{Remark}
\newcommand{\R}{\mathbb R}
\newcommand{\D}{\mathcal D}
\newcommand{\dd}{\,\mathrm d}
\DeclareMathOperator{\conv}{conv}
\DeclareMathOperator{\Var}{Var}
\DeclareMathOperator{\tr}{tr}
\title[The Laugesen--Siudeja eigenvalue conjecture]
{On the Laugesen--Siudeja Conjecture: A Sharp Lower Bound and Failure of the Upper Bound}
\author[Zhengchao Ji]{%
Zhengchao Ji%
\texorpdfstring{%
\protect\\[5pt]
\protect\normalfont\protect\small
Department of Mathematics, China Jiliang University%
\protect\\
\protect\normalfont\protect\small
Hangzhou 310027, People's Republic of China%
\protect\\
\protect\normalfont\protect\small
Correspondence: jizhengchao@zju.edu.cn%
}{}%
}


\begin{document}
\begin{abstract}
We prove the sharp lower bound in the Laugesen--Siudeja conjecture and
give a counterexample to its upper bound. For every bounded convex
planar domain, the first Dirichlet eigenvalue $\lambda_1$, the area $A$,
and the centroidal polar moment of inertia $I$ satisfy
$\lambda_1A^3/I>9\pi^2/2$. The constant is approached by thin isosceles
triangles. The lower bound follows from a rectangle comparison in the
principal directions of the covariance matrix and a planar random-triangle
inequality. For the upper bound, we exhibit an isosceles trapezoid on
which $\lambda_1A^3/I>12\pi^2$, using Temple's inequality and exact
arithmetic.
\end{abstract}
\maketitle

\section{Introduction and main results}

For a bounded connected Lipschitz domain $\Omega\subset\R^2$, write
\[
 0<\lambda_1(\Omega)<\lambda_2(\Omega)\leq\lambda_3(\Omega)\leq\cdots
\]
for the eigenvalues of the positive Dirichlet Laplacian, counted with
multiplicity. We denote the area, centroid, and centroidal polar moment
of inertia by
\begin{equation}\label{eq:definitions}
 A(\Omega)=|\Omega|,\qquad
 c(\Omega)=\frac1{A(\Omega)}\int_\Omega x\dd x,\qquad
 I(\Omega)=\int_\Omega |x-c(\Omega)|^2\dd x
\end{equation}
The quantity
\[
 Q(\Omega)=\lambda_1(\Omega)\frac{A(\Omega)^3}{I(\Omega)}
\]
is invariant under similarities.

P\'olya and Szeg\H{o}'s inclusion method gives upper and lower bounds
for $Q$ on convex domains~\cite[Sections~1.19 and~5.11b]{PSzego}.
Laugesen and Siudeja used rotational symmetry and tight frames to show
that $(\lambda_1+\cdots+\lambda_n)A^3/I$ is maximized by the equilateral
triangle among triangles, by the square among parallelograms, and by
the disk among ellipses~\cite{LS}. These results extend P\'olya's
bounds for the first eigenvalue. When $n=1$, all rectangles and the
equilateral triangle give the same value, $12\pi^2$. This led Laugesen
and Siudeja to formulate the following conjecture.

\begin{conjecture}[Laugesen--Siudeja, {\cite[Conjecture~4.1]{LS}}]
For every bounded convex planar domain $\Omega$,
\begin{equation}\label{eq:conjecture}
 \frac92\pi^2<Q(\Omega)\leq12\pi^2.
\end{equation}
Equality on the right holds for equilateral triangles and all rectangles.
The lower endpoint is approached by degenerate acute isosceles triangles
and sectors.
\end{conjecture}

Freitas obtained more precise eigenvalue asymptotics for thin
triangles~\cite{FreitasThin}. The conjecture was restated in
Laugesen's survey~\cite[Conjecture~7.1]{LaugesenSurvey} and in the
chapter by Laugesen and Siudeja in \emph{Shape Optimization and
Spectral Theory}~\cite[Conjecture~6.62]{LSChapter}.

Two related results use different geometric normalizations. For
quadrilaterals with side lengths $\ell_1,\ldots,\ell_4$, Freitas and
Siudeja~\cite[Theorem~1.1]{FS} proved the sharp bound
$\lambda_1 A^2/(\ell_1^2+\cdots+\ell_4^2)\leq\pi^2/2$.
The spectral bounds for starlike domains in \cite{LSStarlike}
involve inertia and additional geometric factors.

The lower inequality in \eqref{eq:conjecture} holds with its conjectured
constant. The upper inequality fails for the following trapezoid.

\begin{theorem}\label{thm:main}
Consider the convex isosceles trapezoid
\begin{equation}\label{eq:domain}
 \Omega=\left\{(x,y)\in\R^2:
 -\frac75<y<\frac75,\quad |x|<1+\frac{5y}{14}\right\}.
\end{equation}
Then $Q(\Omega)>12\pi^2$. More precisely,
\begin{equation}\label{eq:mainbounds}
 \lambda_1(\Omega)>\frac{77}{20},\qquad
 \frac{12\pi^2 I(\Omega)}{A(\Omega)^3}
 =\frac{457\pi^2}{1176}<\frac{19177}{5000}<\frac{77}{20}.
\end{equation}
\end{theorem}

For the lower endpoint, let $\mathcal C$ be the class of nonempty
bounded convex open subsets of $\R^2$.

\begin{theorem}\label{thm:general-lower}
Every $\Omega\in\mathcal C$ satisfies
\begin{equation}\label{eq:general-lower}
 \lambda_1(\Omega)\frac{A(\Omega)^3}{I(\Omega)}>\frac92\pi^2.
\end{equation}
Moreover,
\[
 \inf_{\Omega\in\mathcal C}
 \lambda_1(\Omega)\frac{A(\Omega)^3}{I(\Omega)}=\frac92\pi^2.
\]
The infimum is unattained. The values on acute isosceles triangles
converge to it as the triangles collapse to a segment.
\end{theorem}

\subsection{Proof strategy}

For the trapezoid, we compute the Rayleigh quotient and the normalized
squared residual of an explicit trial function whose pullback to the
square is polynomial. The spectral moments are rational linear
combinations of $1$ and $\log 3$. Rational bounds for $\log 3$ and a
lower bound for $\lambda_2$ yield, through Temple's inequality, a
lower bound for $\lambda_1$ above the proposed threshold.

For the lower bound, write $C$ for the covariance matrix of the uniform
distribution on $\Omega$, so that $I=A\tr C$. The projected densities
along the principal directions are concave. A lower bound for their
variances controls the widths of the supporting rectangle, whose
first eigenvalue gives
\[
 Q(\Omega)\geq\frac{\pi^2A^2}{24\det C}.
\]
The planar random-simplex inequality of Campi, Colesanti, and
Gronchi~\cite{CCG} gives $\det C/A^2\leq1/108$. Saroglou~\cite{Saroglou}
established that the maximizers are precisely the triangles. Domain
monotonicity gives the strict spectral bound for domains that differ
from their supporting rectangles. Rectangles give $Q=12\pi^2$
directly. A comparison with two rectangles gives the limiting value
for thin triangles. For a convex quadrilateral, the determinant can
be computed from the ratios in which its diagonals divide one
another. This yields the bound in Theorem~\ref{thm:quad-lower}.

\section{The trapezoid and its second eigenvalue}\label{sec:geometry}

Let $\Omega$ be the trapezoid in \eqref{eq:domain}. The vertices of its
closure, in counterclockwise order, are
\[
 \left(-\frac12,-\frac75\right),\quad
 \left(\frac12,-\frac75\right),\quad
 \left(\frac32,\frac75\right),\quad
 \left(-\frac32,\frac75\right).
\]
We write $A$, $c$, and $I$ for its area, centroid, and centroidal polar
moment of inertia.

\begin{lemma}\label{lem:geometry}
The area, centroid, and moment of inertia of $\Omega$ are
\[
 A=\frac{28}{5},\qquad c=\left(0,\frac7{30}\right),\qquad
 I=\frac{6398}{1125},\qquad \frac{I}{A^3}=\frac{457}{14112}.
\]
\end{lemma}
\begin{proof}
Set $h=7/5$. The horizontal section at height $y$ is centered at
$x=0$ and has half-width $1+y/(2h)$. Integrating over these sections gives
\[
 A=4h=\frac{28}{5},\qquad
 \int_\Omega x\dd x\dd y=0,\qquad
 \int_\Omega y\dd x\dd y=\frac{2h^2}{3}.
\]
Thus $c_y=h/6=7/30$. Similarly,
\[
 \int_\Omega x^2\dd x\dd y=\frac{5h}{3}=\frac73,
 \qquad
 \int_\Omega y^2\dd x\dd y=\frac{4h^3}{3}=\frac{1372}{375}.
\]
The formula for $I$ follows by subtracting $A|c|^2=343/1125$ from
the sum of the two second moments. Dividing by $A^3$ gives the
last identity.
\end{proof}

\begin{lemma}\label{lem:second}
For the domain \eqref{eq:domain},
\begin{equation}\label{eq:second}
 \lambda_2(\Omega)>\frac{323}{50}.
\end{equation}
\end{lemma}
\begin{proof}
By the Hong--Krahn--Szeg\H{o} inequality~\cite{Krahn1926},
\begin{equation}\label{eq:hks}
 \lambda_2(\Omega)\geq\frac{2\pi j_{0,1}^2}{A},
\end{equation}
where $j_{0,1}$ is the first positive zero of $J_0$.
Indeed, on each nodal domain of a second eigenfunction, the first
Dirichlet eigenvalue equals $\lambda_2$. The Faber--Krahn
inequality~\cite{Faber1923,Krahn1925} gives an area of at least
$\pi j_{0,1}^2/\lambda_2$ for each such domain. There are nodal domains
of both signs, and their combined area is at most $A$.

We show that $j_{0,1}>12/5$. For $0<x\leq12/5$, the successive
absolute terms of the series for $J_1(x)$ have ratio at most $18/25$.
The alternating-series bound therefore gives
\[
 J_1(x)\geq\frac{x}{2}-\frac{x^3}{16}>0.
\]
Thus $J_0'=-J_1<0$ on $(0,12/5]$. Truncating the series for $J_0$
at $x=12/5$ gives
\begin{equation}\label{eq:besselcheck}
 J_0(12/5)>
 \sum_{k=0}^{7}\frac{(-1)^k(36/25)^k}{(k!)^2}>0.
\end{equation}
It follows that $J_0$ is positive throughout $[0,12/5]$, so
$j_{0,1}>12/5$. Using Lemma~\ref{lem:geometry} in \eqref{eq:hks}, we obtain
\[
 \lambda_2(\Omega)>\frac{72\pi}{35}>
 \frac{72}{35}\frac{314159}{100000}>
 \frac{323}{50}.
\]
For the enclosure $314159/100000<\pi<3927/1250$, use Machin's
identity $\pi=16\arctan(1/5)-4\arctan(1/239)$. Retain the first
$40$ terms of each arctangent series and bound its remainder by the
first omitted term.
\end{proof}

\section{Temple's inequality}\label{sec:temple}

The Rayleigh quotient bounds the first eigenvalue from above.
Temple's inequality~\cite{Temple} gives a lower bound from the same
trial function once its residual and a lower bound for the second
eigenvalue are available.

\begin{lemma}\label{lem:temple}
Let $H$ be a self-adjoint operator, bounded below and with compact
resolvent. Suppose its lowest eigenvalue $\mu_1$ is simple and its
second eigenvalue is $\mu_2$. For a nonzero $v\in\D(H)$, set
\[
 R=\frac{\langle Hv,v\rangle}{\|v\|^2},\qquad
 W=\frac{\|Hv\|^2}{\|v\|^2},\qquad V=W-R^2.
\]
If $R<L\leq\mu_2$, then
\begin{equation}\label{eq:temple}
 \mu_1\geq R-\frac{V}{L-R}.
\end{equation}
\end{lemma}
\begin{proof}
Expand $v$ in an orthonormal eigenbasis. Since
$(\mu_j-\mu_1)(\mu_j-L)\geq0$ for every $j$, summation gives
\[
 0\leq W-(\mu_1+L)R+\mu_1L.
\]
The series converges since $v\in\D(H)$ and hence $\|Hv\|^2<\infty$.
Rearranging the inequality and dividing by $L-R>0$ gives
\eqref{eq:temple}.
\end{proof}

\section{The trial function and its spectral moments}\label{sec:trial}

Let $D=(-1,1)^2$ and define
\begin{equation}\label{eq:map}
 F(s,t)=\left(\frac{2+t}{2}s,\frac75t\right),\qquad
 \det DF=\frac7{10}(2+t).
\end{equation}
The map $F$ takes $D$ onto $\Omega$. Since $2+t$ stays positive near
$\overline D$, the map and its inverse extend smoothly to
neighborhoods of the closures. Define $u$ by
\begin{equation}\label{eq:trial}
 u(F(s,t))=U(s,t):=(1-s^2)(1-t^2)\bigl(p(t)+s^2q(t)\bigr),
\end{equation}
where
\begin{equation}\label{eq:polynomials}
\begin{aligned}
 p(t)={}&81197+44604t-25362t^2-25020t^3
             +5285t^4+4032t^5,\\
 q(t)={}&-19851-33966t-19674t^2+1440t^3
             +8085t^4-10962t^5.
\end{aligned}
\end{equation}
On $\partial D$, one of the factors $1-s^2$ and $1-t^2$ vanishes.
Thus $u$ has zero boundary trace on $\Omega$. It is smooth up to
the boundary and is nonzero since $U(0,0)=81197$. In particular,
$u\in H^2(\Omega)\cap H^1_0(\Omega)\subset\D(-\Delta_D)$.

Define the polynomial
\begin{equation}\label{eq:Hpoly}
 \mathcal H=(784+100s^2)U_{ss}
 -200s(2+t)U_{st}
 +100(2+t)^2U_{tt}+200sU_s.
\end{equation}
In the $(s,t)$ coordinates, differentiation in $x$ and $y$ takes the form
\[
 \partial_x=\frac{2}{2+t}\partial_s,\qquad
 \partial_y=\frac57\left(\partial_t-\frac{s}{2+t}\partial_s\right).
\]
Applying these operators twice yields
\begin{equation}\label{eq:laplace}
 (\Delta u)\circ F=\frac{\mathcal H}{196(2+t)^2}.
\end{equation}
After removing the common factor $7/10$ from the Jacobian, the
integrals for the spectral moments are
\begin{equation}\label{eq:integrals}
\begin{aligned}
 M_0&=\int_D U^2(2+t)\dd s\dd t,\\
 E_0&=-\frac1{196}\int_D\frac{U\mathcal H}{2+t}\dd s\dd t,\\
 W_0&=\frac1{196^2}\int_D\frac{\mathcal H^2}{(2+t)^3}\dd s\dd t.
\end{aligned}
\end{equation}
The change of variables and integration by parts give
\begin{equation}\label{eq:moments}
 R:=\frac{\int_\Omega|\nabla u|^2}{\int_\Omega u^2}
   =\frac{E_0}{M_0},\qquad
 W:=\frac{\int_\Omega|\Delta u|^2}{\int_\Omega u^2}
   =\frac{W_0}{M_0}.
\end{equation}

\begin{proposition}\label{prop:moments}
The normalization integral is
\begin{equation}\label{eq:Mexact}
 M_0=\frac{21810721748942848}{1576575},
\end{equation}
and the moments in \eqref{eq:moments} are
\begin{equation}\label{eq:Rexact}
 R=-\frac{14680948255457463997}{25048250758551552}
 +\frac{2562149793229053885}{4771095382581248}\log3
\end{equation}
and
\begin{equation}\label{eq:Wexact}
 W=-\frac{138049420108616423861415}{58445918436620288}
 +\frac{71805025093027032663285}{33397667678068736}\log3.
\end{equation}
Moreover, with $V=W-R^2$,
\begin{equation}\label{eq:certifiedbounds}
 \frac{77291}{20000}<R<\frac{48307}{12500},\qquad
 0<V<\frac{1871}{50000}.
\end{equation}
\end{proposition}
\begin{proof}
The integrands in \eqref{eq:integrals} can be expanded in monomials
in $s$ and $t$. Integrate first in $s$, using
\[
 \int_{-1}^{1}s^i\dd s=
 \begin{cases}0,&i\text{ odd},\\2/(i+1),&i\text{ even}.
 \end{cases}
\]
For the remaining integrals, set $z=2+t$ and use
\begin{equation}\label{eq:monomialintegration}
 \int_{-1}^{1}\frac{t^j}{(2+t)^m}\dd t
 =\sum_{k=0}^{j}\binom jk(-2)^{j-k}B_{k-m},
 \quad
 B_\ell=\begin{cases}
 \log3,&\ell=-1,\\
 (3^{\ell+1}-1)/(\ell+1),&\ell\ne-1.
 \end{cases}
\end{equation}
Integrating the polynomial $U^2(2+t)$ term by term gives
\eqref{eq:Mexact}. In \eqref{eq:monomialintegration}, take
$m=1$ for $E_0$ and $m=3$ for $W_0$. Evaluating the resulting sums
and dividing by $M_0$ gives \eqref{eq:Rexact} and \eqref{eq:Wexact}.

For the numerical bounds, use the partial sum and remainder bound
\begin{equation}\label{eq:logbounds}
 S_N=\sum_{k=0}^{N-1}\frac1{4^k(2k+1)},\qquad
 \varepsilon_N=\frac4{3(2N+1)4^N}.
\end{equation}
Expanding $\log3=2\operatorname{arctanh}(1/2)$ and bounding the tail
by a geometric series gives $S_N<\log3<S_N+\varepsilon_N$.
Write \eqref{eq:Rexact} as $a+b\log3$ and
\eqref{eq:Wexact} as $c+d\log3$, so that $b,d>0$. For $N=40$,
set
\[
 r_-=a+bS_{40},\quad r_+=a+b(S_{40}+\varepsilon_{40}),\qquad
 w_-=c+dS_{40},\quad w_+=c+d(S_{40}+\varepsilon_{40}).
\]
Direct comparison of these rational numbers gives
\begin{equation}\label{eq:finitechecks}
 \frac{77291}{20000}<r_-\leq r_+<\frac{48307}{12500},\qquad
 0<w_--r_+^2\leq w_+-r_-^2<\frac{1871}{50000}.
\end{equation}
Because $0<r_-<R<r_+$ and $w_-<W<w_+$, the two bounds for
$W-R^2$ in \eqref{eq:finitechecks} imply \eqref{eq:certifiedbounds}.
\end{proof}

\section{Proof of the counterexample}\label{sec:counterexample}

\begin{proof}[Proof of Theorem~\ref{thm:main}]
Apply Lemma~\ref{lem:temple} to the Dirichlet Laplacian and the trial
function in \eqref{eq:trial}. Lemma~\ref{lem:second} allows us to take
$L=323/50$, and Proposition~\ref{prop:moments} gives $R<L$. Hence
\begin{align*}
 \lambda_1(\Omega)
 &\geq R-\frac{V}{L-R}\\
 &>\frac{77291}{20000}
 -\frac{1871/50000}{323/50-48307/12500}
 >\frac{77}{20}.
\end{align*}
The last step is the rational comparison
\[
 \frac{77291}{20000}-\frac{77}{20}
 >\frac{1871/50000}{323/50-48307/12500}.
\]
The geometric threshold is estimated using Lemma~\ref{lem:geometry}
and the upper bound for $\pi$
\[
 \frac{12\pi^2I}{A^3}
 =\frac{457\pi^2}{1176}
 <\frac{457}{1176}\left(\frac{3927}{1250}\right)^2
 <\frac{19177}{5000}<\frac{77}{20}.
\]
This proves \eqref{eq:mainbounds} and gives
$\lambda_1(\Omega)>12\pi^2 I/A^3$.
\end{proof}

\begin{remark}
In terms of $Q$, the estimate $\lambda_1>77/20$ gives
\[
 \lambda_1(\Omega)\frac{A^3}{I}>
 \frac{77}{20}\frac{14112}{457}=\frac{271656}{2285}.
\]
\end{remark}

\section{The sharp lower bound for convex planar domains}
\label{sec:general-lower}

To compare $\Omega$ with its supporting rectangle, we need a lower
bound for the variance of a concave density in terms of the length
of its interval.

\begin{lemma}\label{lem:quad-variance}
If $f$ is a nonnegative continuous concave probability density on an
interval of length $L>0$, then its variance is at least $L^2/24$.
The constant $1/24$ is optimal in this class.
\end{lemma}
\begin{proof}
By rescaling we may take the interval to be $[0,1]$. Assume first
that $f$ is piecewise affine. For $0<t<1$, set
\[
 h_t(x)=\begin{cases}2x/t,&0\leq x\leq t,\\
 2(1-x)/(1-t),&t\leq x\leq1.
 \end{cases}
\]
Also let $h_0(x)=2(1-x)$ and $h_1(x)=2x$.
Direct integration gives, for all $0\leq t\leq1$,
\[
 \int_0^1 h_t=1,\quad
 \int_0^1 xh_t(x)\,dx=\frac{1+t}{3},\quad
 \int_0^1x^2h_t(x)\,dx=\frac{1+t+t^2}{6}.
\]
Consequently
\[
 \Var(h_t)=\frac{1-t+t^2}{18}
 =\frac1{24}+\frac{(t-1/2)^2}{18}\geq\frac1{24}.
\]
Let $t_j$ run over the interior breakpoints of $f$, and let
$\kappa_j\geq0$ be the decrease in slope at $t_j$. The function
\[
 G(x,t)=\min(x,t)\bigl(1-\max(x,t)\bigr)
\]
vanishes at both endpoints and has a unit downward jump in slope at
$t$. The endpoint values and slope jumps determine $f$, so
\[
 f(x)=f(0)(1-x)+f(1)x+\sum_j\kappa_jG(x,t_j).
\]
Using $G(x,t)=t(1-t)h_t(x)/2$, we can write
\[
 f(x)=\frac{f(0)}2h_0(x)+\frac{f(1)}2h_1(x)
     +\sum_j\frac{\kappa_jt_j(1-t_j)}2h_{t_j}(x).
\]
The coefficients are nonnegative and sum to one, since $f$ and each
$h_t$ have integral one. By the law of total variance, the variance
of this mixture is at least the weighted average of the component
variances. Each of these is at least $1/24$.

For a continuous concave density, take piecewise linear interpolants
on successively finer uniform partitions of $[0,1]$ and normalize
them to have integral one. They are nonnegative and concave, and
uniform convergence gives convergence of their first two moments.
The variance bound passes to the limit and rescales to $L^2/24$.
The density $h_{1/2}$ attains the constant on $[0,1]$.
\end{proof}

Let $\Omega\in\mathcal C$. Write $A=|\Omega|$ and define its
covariance matrix by
\[
 C=\frac1A\int_\Omega(x-c(\Omega))(x-c(\Omega))^T\dd x.
\]
The matrix $C$ is positive definite and $I(\Omega)=A\tr C$.
The closure $K=\overline\Omega$ is a convex body with interior $\Omega$.
Since its boundary has measure zero, we may use either set to compute
the area, centroid, and covariance.

\begin{lemma}\label{lem:general-rectangle}
Let $s_1,s_2>0$ be the eigenvalues of $C$, and let $W_1,W_2$ be the
projection widths of $K$ along corresponding orthonormal eigenvectors.
If $R$ is the open supporting rectangle in these directions, then
\begin{equation}\label{eq:general-rectangle}
 \lambda_1(\Omega)\geq\lambda_1(R)
 \geq\frac{\pi^2\tr C}{24\det C}.
\end{equation}
Consequently,
\begin{equation}\label{eq:general-covariance-spectral}
 \lambda_1(\Omega)\frac{A^3}{I(\Omega)}
 \geq\frac{\pi^2A^2}{24\det C}.
\end{equation}
\end{lemma}
\begin{proof}
Fix a projection direction and let $J(t)$ be the perpendicular section
at coordinate $t$. The projected density is $|J(t)|/A$. Convexity gives
\[
 (1-\theta)J(t)+\theta J(u)
 \subset J\bigl((1-\theta)t+\theta u\bigr),\qquad 0\leq\theta\leq1.
\]
Minkowski addition of intervals adds their lengths, so the inclusion
implies concavity of the projected density on the interior of the
projection interval. Boundedness gives finite one-sided limits
at the endpoints, which define a continuous concave extension.
The values at the two endpoints do not affect the variance.

Applying Lemma~\ref{lem:quad-variance} in the two principal
directions gives $s_i\geq W_i^2/24$ for $i=1,2$.

Since $\Omega\subset R$, domain monotonicity and the rectangle spectrum
give
\[
 \lambda_1(\Omega)\geq\lambda_1(R)
 =\pi^2\left(W_1^{-2}+W_2^{-2}\right)
 \geq\frac{\pi^2}{24}\left(s_1^{-1}+s_2^{-1}\right)
 =\frac{\pi^2\tr C}{24\det C}.
\]
Multiplication by $A^2/\tr C$ proves the last assertion.
\end{proof}

The geometric input for the determinant bound is the maximality
of triangles for the normalized second moment of random triangle
area~\cite{CCG,Saroglou}.

\begin{proposition}
\label{prop:general-determinant}
For every planar convex body $K$ of area $A$ and covariance matrix $C$,
\begin{equation}\label{eq:general-determinant}
 \frac{\det C}{A^2}\leq\frac1{108}.
\end{equation}
\end{proposition}
\begin{proof}
For independent uniform points $X_0,X_1,X_2$ in $K$, put
\[
 \mathcal S(K)=\frac1{A^2}\,
 \mathbb E\!\left[\operatorname{Area}
       (\conv\{X_0,X_1,X_2\})^2\right].
\]
In the notation of Campi, Colesanti, and Gronchi, this is
$M_2(K;3)$ in dimension two. The factor $A^{-2}$ in the expectation
corresponds to $A^{-5}$ in the integral over $K^3$.
Their planar inequality~\cite[Section~3, inequality~(11)]{CCG} gives
\begin{equation}\label{eq:random-triangle-input}
 \mathcal S(K)\leq\mathcal S(T)
 \quad\text{for every nondegenerate triangle }T.
\end{equation}
They prove the inequality for polygons and extend it to planar
convex bodies by continuity.

To express $\mathcal S$ in terms of $C$, center the points by setting
$Y_i=X_i-\mathbb E X_i$, and write $[x,y]=\det[x,y]$. Twice the signed
triangle area is
\[
 D=[Y_1,Y_2]+[Y_2,Y_0]+[Y_0,Y_1].
\]
Each summand has mean zero and second moment
\[
 \mathbb E\!\left([Y_i,Y_j]^2\right)
 =2(C_{11}C_{22}-C_{12}^2)=2\det C\qquad(i\ne j).
\]
For each mixed product, conditioning on the shared point gives
zero expectation by independence and centering. Therefore
$\mathbb E D^2=6\det C$ and
\begin{equation}\label{eq:random-triangle-covariance}
 \mathcal S(K)=\frac{3\det C}{2A^2}.
\end{equation}

For $T=\conv\{(0,0),(1,0),(0,1)\}$, direct integration gives
\[
 A_T=\frac12,\qquad c_T=(1/3,1/3),\qquad
 C_T=\begin{pmatrix}1/18&-1/36\\-1/36&1/18\end{pmatrix}.
\]
Thus $\det C_T=1/432$, $\det C_T/A_T^2=1/108$, and
$\mathcal S(T)=1/72$. Substitution in
\eqref{eq:random-triangle-input}--\eqref{eq:random-triangle-covariance}
proves \eqref{eq:general-determinant}.
\end{proof}

\begin{proof}[Proof of the lower inequality in Theorem~\ref{thm:general-lower}]
Lemma~\ref{lem:general-rectangle} and
Proposition~\ref{prop:general-determinant} give
\[
 \lambda_1(\Omega)\frac{A^3}{I(\Omega)}
 \geq\frac{\pi^2A^2}{24\det C}\geq\frac92\pi^2.
\]
Suppose that $\Omega$ is a proper subset of the supporting rectangle
$R$. Since $\Omega=\operatorname{int}\overline\Omega$, the open set
$R$ cannot be contained in $\overline\Omega$. Thus
$R\setminus\overline\Omega$ contains an open ball. If
$\lambda_1(\Omega)=\lambda_1(R)$, extending a first eigenfunction on
$\Omega$ by zero would give a minimizer of the Rayleigh quotient
in $H^1_0(R)$. Every first eigenfunction on $R$ is a nonzero
multiple of the product of the two first sine modes and is nowhere
zero in $R$. The open ball on which the extension vanishes rules
out equality. Thus the first inequality in \eqref{eq:general-rectangle}
is strict.

If $\Omega=R$, with side lengths $L,W$, then
\[
 A=LW,\qquad I=\frac{LW(L^2+W^2)}{12},\qquad
 \lambda_1=\pi^2(L^{-2}+W^{-2}),
\]
so $\lambda_1A^3/I=12\pi^2>9\pi^2/2$.
\end{proof}

\subsection{Optimality of the lower constant}
\label{subsec:sharpness}

The original conjecture identifies thin acute isosceles triangles
as a limiting family~\cite[Conjecture~4.1]{LS}. The leading spectral
term can be recovered directly by enclosing each triangle in a
rectangle and placing a second rectangle near its widest section.
For the finer asymptotics, see \cite{FreitasThin}.

\begin{proposition}\label{prop:sharpness}
For $0<\varepsilon<1$, let
\[
 T_\varepsilon=\{(x,y):0<x<1,\ |y|<\varepsilon x\}.
\]
These are acute isosceles triangles, and
\begin{equation}\label{eq:sharpness-limit}
 \lim_{\varepsilon\downarrow0}
 \lambda_1(T_\varepsilon)
 \frac{A(T_\varepsilon)^3}{I(T_\varepsilon)}
 =\frac92\pi^2.
\end{equation}
Consequently, over all nonempty bounded convex planar open sets,
\begin{equation}\label{eq:sharp-infimum}
 \inf_\Omega\lambda_1(\Omega)\frac{A(\Omega)^3}{I(\Omega)}
 =\frac92\pi^2,
\end{equation}
and the infimum is not attained.
\end{proposition}
\begin{proof}
Integration over vertical sections gives
\[
 A_\varepsilon=\varepsilon,\qquad
 c_\varepsilon=(2/3,0),\qquad
 \frac{I_\varepsilon}{A_\varepsilon}
 =\frac1{18}+\frac{\varepsilon^2}{6}.
\]
Indeed, the averages of $x,x^2,y^2$ are respectively
$2/3,1/2,\varepsilon^2/6$. In particular,
\begin{equation}\label{eq:thin-triangle-normalization}
 \lambda_1(T_\varepsilon)\frac{A_\varepsilon^3}{I_\varepsilon}
 =\frac{18\varepsilon^2\lambda_1(T_\varepsilon)}
        {1+3\varepsilon^2}.
\end{equation}

For $0<\delta<1$, the rectangles
\[
 R_{\varepsilon,\delta}
 =(1-\delta,1)\times
  \bigl(-\varepsilon(1-\delta),\varepsilon(1-\delta)\bigr),
 \qquad
 S_\varepsilon=(0,1)\times(-\varepsilon,\varepsilon)
\]
satisfy $R_{\varepsilon,\delta}\subset T_\varepsilon\subset S_\varepsilon$.
Thus the exact rectangle spectrum and domain monotonicity yield
\[
 \pi^2\left(1+\frac1{4\varepsilon^2}\right)
 \leq\lambda_1(T_\varepsilon)
 \leq\pi^2\left(\frac1{\delta^2}
                +\frac1{4\varepsilon^2(1-\delta)^2}\right).
\]
Take $\delta=\varepsilon^{2/3}$ and multiply by $\varepsilon^2$ to obtain
\begin{equation}\label{eq:thin-triangle-squeeze}
 \pi^2\left(\frac14+\varepsilon^2\right)
 \leq\varepsilon^2\lambda_1(T_\varepsilon)
 \leq\pi^2\left(\varepsilon^{2/3}
         +\frac1{4(1-\varepsilon^{2/3})^2}\right).
\end{equation}
Both bounds tend to $\pi^2/4$ as $\varepsilon\downarrow0$.
Substitution in \eqref{eq:thin-triangle-normalization} gives
\eqref{eq:sharpness-limit}. Together with the strict bound
\eqref{eq:general-lower}, this proves \eqref{eq:sharp-infimum}
and shows that the infimum is unattained.
\end{proof}

\section{A quantitative bound for convex quadrilaterals}\label{sec:quad-lower}

For a quadrilateral, the ratios in which the diagonals divide one
another determine $\det C/A^2$. We use this identity to refine the
general bound. In this section $K$ is the interior of a quadrilateral
with four vertices in convex position, $A=|K|$, and $C$ is its
covariance matrix.

\begin{theorem}\label{thm:quad-lower}
Let $a,c>0$ be the lengths of the two portions of one diagonal, and let
$b,d>0$ be the lengths of the two portions of the other diagonal. Define
\[
 u=\frac{c-a}{a+c},\qquad v=\frac{b-d}{b+d},\qquad
 \Delta=(1-u^2)(1-v^2).
\]
Then $0<\Delta\leq1$ and
\[
 \frac{\lambda_1(K)A^3}{I(K)}
 \geq\frac{18\pi^2}{4-\Delta}>\frac92\pi^2.
\]
\end{theorem}

\begin{lemma}\label{lem:quad-det}
For the quadrilateral in Theorem~\ref{thm:quad-lower},
\[
 \frac{\det C}{A^2}=\frac{4-\Delta}{432}.
\]
\end{lemma}
\begin{proof}
Translate the intersection of the diagonals to the origin. Choose
unit vectors $p,q$ along them so that the vertices are $-ap,-dq,cp,bq$
in counterclockwise order and $\det[p,q]>0$. Set $B=[p,q]$. In the
coordinates $(\xi,\eta)$ associated with $B$, the vertices are
$(-a,0),(0,-d),(c,0),(0,b)$. The area in these coordinates is
$A_0=(a+c)(b+d)/2$, so $A=A_0\det B$. The uniform distribution has
moments
\[
 \begin{gathered}
 \mathbb E\xi=\frac{c-a}{3},\qquad
 \mathbb E\eta=\frac{b-d}{3},\\
 \mathbb E\xi^2=\frac{a^2-ac+c^2}{6},\qquad
 \mathbb E\eta^2=\frac{b^2-bd+d^2}{6},\\
 \mathbb E\xi\eta=\frac{(c-a)(b-d)}{12}.
 \end{gathered}
\]
To obtain these formulas, split the quadrilateral into four triangles
at the origin. For vertices $(0,0),(r,0),(0,s)$, the averages of
$\xi,\eta,\xi^2,\eta^2,\xi\eta$ are $r/3,s/3,r^2/6,s^2/6,rs/12$,
respectively. Average over $r\in\{-a,c\}$ and $s\in\{-d,b\}$ with
area weights $|rs|/2$, then subtract the products of the means to get
\[
 C_0=\begin{pmatrix}
 (a^2+ac+c^2)/18&-(c-a)(b-d)/36\\
 -(c-a)(b-d)/36&(b^2+bd+d^2)/18
 \end{pmatrix}.
\]
Since $C=BC_0B^T$, we have $\det C/A^2=\det C_0/A_0^2$.
Set $X=a+c$ and $Y=b+d$. Then
\[
 C_0=\begin{pmatrix}
 X^2(3+u^2)/72&-XYuv/36\\
 -XYuv/36&Y^2(3+v^2)/72
 \end{pmatrix}.
\]
Taking the determinant and dividing by $A_0^2=X^2Y^2/4$ gives
\[
 \frac{\det C}{A^2}
 =\frac{(3+u^2)(3+v^2)-4u^2v^2}{1296}
 =\frac{3+u^2+v^2-u^2v^2}{432}
 =\frac{4-(1-u^2)(1-v^2)}{432}.
\]
\end{proof}

\begin{proof}[Proof of Theorem~\ref{thm:quad-lower}]
By Lemma~\ref{lem:general-rectangle} and Lemma~\ref{lem:quad-det},
\[
 \frac{\lambda_1(K)A^3}{I(K)}
 \geq\frac{\pi^2}{24}\frac{432}{4-\Delta}
 =\frac{18\pi^2}{4-\Delta}.
\]
The diagonals meet in the interior of both segments, so $|u|,|v|<1$.
Consequently $0<\Delta\leq1$, and the bound is strictly larger than
$9\pi^2/2$.
\end{proof}

The estimate gives the explicit gap
\[
 \frac{18\pi^2}{4-\Delta}-\frac92\pi^2
 =\frac{9\pi^2\Delta}{2(4-\Delta)}.
\]

\section{Open problem}

Theorem~\ref{thm:general-lower} identifies the lower endpoint of $Q$
on $\mathcal C$. The infimum is $9\pi^2/2$, approached by thin acute
isosceles triangles. At the upper end, the bounds of Laugesen and
Siudeja~\cite{LS} give the maximum $12\pi^2$ within each of the classes
of triangles and parallelograms. The equilateral triangle and every
rectangle attain this value.

For the isosceles trapezoid in \eqref{eq:domain}, the spectral and
geometric estimates in \eqref{eq:mainbounds} give
$Q(\Omega)>271656/2285$. Combining this lower bound for the supremum
with the inclusion bounds of P\'olya and Szeg\H{o}
\cite[Sections~1.19 and~5.11b]{PSzego}, we have
\[
 12\pi^2<\frac{271656}{2285}
 <\sup_{\Omega\in\mathcal C}Q(\Omega)<\infty.
\]
In particular, any domain attaining this supremum lies outside both
of the classes above.

\begin{openproblem}
Determine the sharp upper bound for $Q$ on $\mathcal C$. Is it attained?
If so, characterize the maximizing domains up to similarity.
\end{openproblem}

\section*{Declaration of competing interests}
The author declares no competing interests.

\section*{Data availability}
No datasets were used in this study.

\section*{Funding}
This work was supported by the National Natural Science Foundation of
China (grant number 12501068) and by the Natural Science Foundation of
Zhejiang Province (grant number LQN25A010001).

\section*{Acknowledgements}
The author is grateful to Professors Kefeng Liu and Hongwei Xu for their
continuous support, valuable advice, and encouragement.

\noindent\textit{AI assistance.}
All proof strategies in this article were developed by the author.
The author used OpenAI Codex (GPT 6) to assist with numerical simulations,
the construction of the counterexample, LaTeX typesetting, and language
editing. The author independently reviewed all AI-assisted work
and takes full responsibility for the final manuscript.

\end{document}